\RequirePackage{fix-cm}

\documentclass[a4paper]{article}
\usepackage{graphicx}
\usepackage{amsmath}
\usepackage{amssymb}
\usepackage{float}

\usepackage{subfig}
\DeclareMathOperator*{\argmin}{arg\,min}
\newcommand*{\QED}{\hfill\ensuremath{\square}}%
\newcommand*{\QEDM}{\tag*{$\square$}}%
\newtheorem{remark}{Remark}
\newtheorem{example}{Example}
\newtheorem{theorem}{Theorem}
\newtheorem{lemma}{Lemma}

\newtheorem{property}{Property}
\newtheorem{proposition}{Proposition}

\providecommand{\keywords}[1]
{
  \small	
  \textbf{\textit{Keywords---}} #1
}

\newcommand{\arrivalrate}{\lambda}
\newcommand{\servicereq}{B}

\newcommand{\servicereqr}{b}

\newcommand{\latency}{T}
\newcommand{\workload}{V}
\newcommand{\workloadstate}{\boldsymbol{\omega}}
\newcommand{\workloadserver}{\omega}
\newcommand{\workloadstateaux}{\boldsymbol{\tilde{\omega}}}
\newcommand{\workloadserveraux}{\tilde{\omega}}
\newcommand{\waitingtime}{W}

\newcommand{\nrservers}{N}

\newcommand{\nrreplicas}{d}

\begin{document}

\title{Redundancy scheduling with scaled Bernoulli service requirements
}



\author{Youri Raaijmakers$^{1,}$\footnote{Corresponding author: Youri Raaijmakers (y.raaijmakers@tue.nl)}, Sem Borst$^{1,2}$, Onno Boxma$^{1}$ \\ \\
$^{1}$\small Department of Mathematics and Computer Science,\\ \small Eindhoven University of Technology, The Netherlands \\
$^{2}$\small Nokia Bell Labs, Murray Hill, USA
}




\maketitle

\begin{abstract}
Redundancy scheduling has emerged as a powerful strategy for improving
response times in parallel-server systems.
The key feature in redundancy scheduling is replication of a job
upon arrival by dispatching replicas to different servers.
Redundant copies are abandoned as soon as the first of these replicas
finishes service.
By creating multiple service opportunities, redundancy scheduling
increases the chance of a fast response from a server that is quick
to provide service, and mitigates the risk of a long delay incurred
when a single selected server turns out to be slow.

The diversity enabled by redundant requests has been found to strongly
improve the response time performance, especially in case of highly
variable service requirements.
Analytical results for redundancy scheduling are unfortunately scarce
however, and even the stability condition has largely remained elusive
so far, except for exponentially distributed service requirements.
In order to gain further insight in the role of the service requirement
distribution, we explore the behavior of redundancy scheduling
for scaled Bernoulli service requirements.
We establish a sufficient stability condition for generally distributed service requirements and we show that, for scaled Bernoulli service requirements, this condition is also asymptotically nearly necessary. 
This stability condition differs drastically
from the exponential case, indicating that the stability condition
depends on the service requirements in a sensitive and intricate manner.
\end{abstract}

\section{Introduction}
\label{sec: introduction}
Redundancy scheduling has recently attracted strong interest as a strategy
for significantly reducing response times in parallel-server systems
\cite{APS-ESM,APS-SMDR,ABV-UFRM,GHBSW-DSSJS,GHBSWVZ-PDR,GZDHBHSW-RLR,GZDHBHSW-QR,PC-CERPS,SLR-WRRRL,SXYZ-TBPPS,VGMSRS-LR}. 
The key feature in redundancy scheduling is replication of a job
upon arrival allowing replicas to be assigned to, say, $\nrreplicas$
different servers, chosen uniformly at random (without replacement).
Redundant replicas are abandoned as soon as the first of these replicas
either starts service (`cancel-on-start' or c.o.s.) or completes service
(`cancel-on-completion' or c.o.c.).
By creating multiple service opportunities, redundancy scheduling boosts
the chance of a fast response from a server that is swift to provide
service, and alleviates the risk of a long delay incurred when a job is
assigned to a single server that may be slow.
Note that the c.o.c.\ and c.o.s.\ policies both ensure that the first
replica starts service at the server with the smallest workload
among the $\nrreplicas$ selected servers.
The possibly concurrent service of multiple replicas under the c.o.c.\
policy provides a further hedge against potentially slow execution
of the first replica in case replicas are independent
(although it may also result in wastage of service effort).

The diversity offered by redundant requests has been shown to strongly
improve the response time performance, especially in case of highly
variable service requirements.
Analytical results for redundancy scheduling are unfortunately scarce
however, and have largely remained limited to exponentially distributed
service requirements.
Specifically, Gardner \textit{et al.}~\cite{GHBSWVZ-PDR} extensively
analyzed the c.o.c.\ redundancy policy with exponentially distributed
service requirements.
They established the stability condition and showed that it does
\textit{not} depend on the number of replicas~$\nrreplicas$,
and thus coincides with the nominal condition without any redundancy.
This may be explained from the fact that even with concurrent service
the expected aggregate amount of time invested in the service of a job
remains equal to the mean service requirement of a single instance
due to the memoryless property of the exponential distribution.
Gardner \textit{et al.}~\cite{GHBSWVZ-PDR} also derived an explicit
expression for the expected latency, and proved that the latency is
decreasing in the number of replicas~$\nrreplicas$.
Simulation experiments additionally demonstrated greater improvements
in the latency in case of highly variable service requirements,
particularly heavy-tailed distributions.

We are not aware of any analytical results for the c.o.c.\ redundancy
policy with independent replicas and nonexponential service requirements.
Hellemans \& Van Houdt~\cite{HH-AIR} consider the c.o.c.\ policy
with \textit{identical} replicas, and derive a differential equation
for the marginal workload distribution at each of the servers
in a limiting regime where the number of servers grows large.
While the differential equation implicitly captures the stability
condition, it does not yield any analytical expression, and the
derivations for identical replicas rely on highly specific arguments
that do not extend to independent replicas.
It is also worth observing that the c.o.s.\ redundancy policy is equivalent
to a power-of-$d$ version of the Join-the-Smallest-Workload policy.
While the workload and waiting-time distributions for these policies
do not appear analytically tractable, the stability condition is simple
and coincides with the nominal condition without any redundancy since
no concurrent service takes place.

In order to gain further insight in the role of the service requirement
distribution, we focus in the present paper on the behavior of the
c.o.c.\ redundancy policy for scaled Bernoulli service requirements.
While this is admittedly a rather special case, it provides a typical
instance of highly variable service requirements for which redundancy
scheduling is particularly relevant, and is also of intrinsic merit
given the paucity of analytical results for general service requirement
distributions.

First of all, we establish a simple sufficient stability condition
in terms of a lower bound for the system capacity, i.e., the maximum aggregate load that can be supported.
The lower bound is obtained from a stochastic coupling between the
maximum workload across all the servers and the workload in a related
single-server queue with the same arrival process and a service
requirement that corresponds to the minimum service requirement
across $d$~replicas.
The lower bound for the system capacity grows
without bound with (a) the `scale' of the service requirement
and (b) the number of replicas~$d$, but remarkably enough (c) does
\textit{not} depend on the number of servers at all
(assuming that number to be at least equal to the number of replicas~$d$).
The `scale' of the service requirement here refers to its non-zero value
relative to its mean and provides a proxy for the degree of variability.
The growth in the system capacity with (a) and (b) reflects the huge
benefits provided by redundancy scheduling for highly variable service
requirements.

In view of~(c), the lower bound may at first sight seem loose
for a larger number of servers, but we will use a further stochastic
comparison argument to prove that it is in fact asymptotically tight
when the scale of the service requirement grows suitably large.
This implies that increasing the number of replicas significantly increases
the system capacity, while adding servers does \textit{not} asymptotically.
Or stated differently, given the number of replicas~$d$,
redundancy scheduling ensures that asymptotically just $d$~servers
suffice to achieve the capacity achievable with any number of servers,
which further highlights the great gains provided by redundancy
scheduling for highly variable service requirements.

The remainder of the paper is organized as follows.
In Section~2 we present a detailed model description, state a sufficient stability condition for generally distributed service requirements.
In Section~3 we prove that this condition is also asymptotically nearly necessary for scaled Bernoulli service requirements.
An upper bound for the expected waiting time is derived
in Section~4 and in Section~5 we provide a conclusion.

\section{Workload model and sufficient stability condition}
\label{sec: workload model}
We consider a system with $\nrservers$ parallel servers. Jobs arrive according to a Poisson process with rate $\arrivalrate$. 
Each arriving job is replicated and immediately allocated to $\nrreplicas$ servers chosen uniformly at random (without replacement). The replicas at each server are served in order of arrival (FCFS) and the job is completed as soon as the first replica finishes service, whereafter the other $\nrreplicas-1$ replicas are instantaneously abandoned. The service requirements of the $\nrreplicas$ replicas are assumed to be independent and identically distributed (i.i.d.) copies of some random variable $\servicereq$. Note that this model corresponds to the independent runtime (IR) model described in \cite{GHBSWVZ-PDR}.

Let $\workloadstate = (\workloadserver_{1},\dots,\workloadserver_{\nrservers})$ denote the workload of the system, where $\workloadserver_{i}$ is the workload at server $i$, for $i=1,\dots,\nrservers$. 
Here we define workload as the \textit{real} amount of work, i.e., the amount of work a server needs to complete to become idle in the absence of any arrivals. This may be smaller than the sum of the service requirements of all the replicas at the server since some replicas may get partly or entirely abandoned, see Example \ref{exam: workloads}.  
Let $s_{j}$ and $\servicereqr_{j}$ denote the sampled server and the realized service requirement of the $j$-th replica, respectively, for $j=1,\dots,\nrreplicas$. The first replica will finish service on server $s_{j^{*}}$, where $j^{*} = \argmin_{j \in \{1,\dots,\nrreplicas \}}(\workloadserver_{s_{j}} + \servicereqr_{j} )$. The workload of server $s_{j}$ is then $\max \{\workloadserver_{s_{j^{*}}} + \servicereqr_{j^{*}}, \workloadserver_{s_{j}}\}$, for $j=1,\dots,\nrreplicas$.
\begin{example}
\label{exam: workloads}
Consider a system with $\nrservers=4$, $\nrreplicas=2$ and workload state\\ $\workloadstate=(4.1,4.1,3.5,2.3)$. Then, after an arrival with service requirements $(2.2,1.5)$ on servers $2$ and $4$, the new workload state is $\workloadstate_{\text{new}}=(4.1,4.1,3.5,3.8)$.
\end{example}

Let $\workloadstate_{(\cdot)}$ denote the workloads arranged in descending order, thus $\workloadstate_{(\cdot)} = \{ \workloadstate \in \mathbb{R}_{+}^{\nrservers}: \workloadserver_{(1)} \geq \workloadserver_{(2)} \geq \hdots \geq \workloadserver_{(\nrservers)} \}$. Throughout this paper we refer to \textit{synchronicity} as the situation in which all workloads are equal, i.e., $\workloadserver_{1}=\hdots=\workloadserver_{\nrservers}$. Moreover, let $\mathcal{S}_{\text{trun}}$ denote the truncated state space of the ordered workload vectors with $\mathcal{S}_{\text{trun}} = \{ \workloadstate \in \mathbb{R}_{+}^{\nrservers}: \workloadserver_{(1)} = \hdots = \workloadserver_{(\nrreplicas)} \geq \workloadserver_{(\nrreplicas+1)} \geq \hdots \geq \workloadserver_{(\nrservers)} \}$.

The next property states that the $\nrreplicas$ largest workloads will always be equal from some point onward. We will later see that under certain conditions the system will in fact be in full synchronicity nearly all the time. 

\begin{property}
\label{prop: synchronicity first d servers}
If $\workloadstate \in \mathcal{S}_{\text{trun}}$ then $\workloadstate_{\text{new}} \in \mathcal{S}_{\text{trun}}$, where $\workloadstate_{\text{new}}$ is any future workload. In other words, once the largest $\nrreplicas$ workloads are equal, they will always remain equal.
\end{property}

\noindent \textbf{Proof:}
Consider the two options, either i) $\min_{j \in \{1,\dots,\nrreplicas\}} (\workloadserver_{s_{j}} + \servicereqr_{j}) \leq \workloadserver_{(1)}=\dots=\workloadserver_{(\nrreplicas)}$ in which case we have $\workloadserver_{\text{new},s_{l}} = \max\{\min_{j \in \{1,\dots,\nrreplicas\}} (\workloadserver_{s_{j}} + \servicereqr_{j}),\workloadserver_{s_{l}}\}\leq \workloadserver_{(1)}$, for $l=1,\dots,\nrreplicas$, therefore $\workloadserver_{(1)}=\dots=\workloadserver_{(\nrreplicas)}=\workloadserver_{\text{new},(1)}=\dots=\workloadserver_{\text{new},(\nrreplicas)}$, ii) $ \min_{j \in \{1,\dots,\nrreplicas\}} (\workloadserver_{s_{j}} + \servicereqr_{j} ) > \workloadserver_{(1)} = \dots = \workloadserver_{(\nrreplicas)}$ in which case $\workloadserver_{\text{new},s_{l}}=\max\{\min_{j \in \{1,\dots,\nrreplicas\}} (\workloadserver_{s_{j}} + \servicereqr_{j}),\workloadserver_{s_{l}}\} = \min_{j \in \{1,\dots,\nrreplicas\}} (\workloadserver_{s_{j}} + \servicereqr_{j})$, for $l=1,\dots,\nrreplicas$, therefore $\workloadserver_{\text{new},s_{1}}=\dots=\workloadserver_{\text{new},s_{\nrreplicas}}=\workloadserver_{\text{new},(1)} = ... = \workloadserver_{\text{new},(\nrreplicas)}$. In both cases $\workloadstate_{\text{new}} \in \mathcal{S}_{\text{trun}}$, thus by a simple induction argument it follows that there are always $\nrreplicas$ servers with the same maximum workload.\\

Before stating and proving a sufficient stability condition, we prove the following lemma for generally distributed service requirements. 

\begin{lemma}
\label{lem: bound maximum workload}
The sequence of maximum workloads $\workloadserver_{(1)}$ at arbitrary epochs is stochastically upper bounded by the sequence of workloads $\workloadserver_{M/G/1}$ in a corresponding $M/G/1$ queue with arrival rate $\arrivalrate_{M/G/1} = \arrivalrate$ and generic service requirement\\ $\servicereq_{M/G/1}=\min \{ \servicereq_{1},\dots,\servicereq_{\nrreplicas}\}$, provided that the initial maximum workload $\workloadserver_{(1)}$ is smaller than the initial workload in the $M/G/1$ queue.
\end{lemma}
\noindent \textbf{Proof:}
The proof follows by induction. Note that for the initial state the statement is satisfied. Assume that $\workloadserver_{(1)} \leq \workloadserver_{M/G/1}$ after the $k$-th arrival. Then, after the $(k+1)$-th arrival the new workload is $\workloadserver_{\text{new},s_{l}} = \max\{\min_{j \in \{1,\dots,\nrreplicas\}} (\workloadserver_{s_{j}} + \servicereqr_{j}),\workloadserver_{s_{l}}\} \leq \max\{\min_{j \in \{1,\dots,\nrreplicas\}} (\workloadserver_{(1)} + \servicereqr_{j}),\workloadserver_{(1)}\} = \workloadserver_{(1)} + \min_{j \in \{1,\dots,\nrreplicas\}} \servicereqr_{j}$, for $l=1,\dots,\nrreplicas$, since $\workloadserver_{i} \leq \workloadserver_{(1)}$ for all $i=1,\dots,\nrservers$. 
Thus the increase in maximum workload is bounded by $\min_{j \in \{1,\dots,\nrreplicas\}} \servicereqr_{j}$, which is exactly the increase in workload in the corresponding $M/G/1$ queue.\QED 

\begin{remark}
\label{rem: workload bound synchronicity}
Observe that in synchronicity, in which all servers have the maximum workload, the bound $\min_{j \in \{1,\dots,\nrreplicas\}} \servicereqr_{j}$ is tight, since here every arrival adds exactly $\min_{j \in \{1,\dots,\nrreplicas\}} \servicereqr_{j}$ work to each of the $\nrreplicas$ sampled servers.
\end{remark}

\begin{proposition}
\label{propo: stability condition IR model}
A sufficient stability condition is
\begin{align}
\arrivalrate \mathbb{E}[\min \{B_{1},\dots,B_{\nrreplicas} \}] < 1.
\label{eq: stability condition proof}
\end{align}
\end{proposition}
\noindent \textbf{Proof:}
By Lemma~\ref{lem: bound maximum workload} we know that the maximum workload in the system is bounded by the workload in a corresponding $M/G/1$ queue with arrival rate $\arrivalrate_{M/G/1} = \arrivalrate$ and generic service requirement $\servicereq_{M/G/1} = \min \{ \servicereq_{1},\dots,\servicereq_{\nrreplicas}\}$.
The (necessary and sufficient) stability condition for the latter $M/G/1$ queue is given by
\begin{align*}
\rho &= \arrivalrate_{M/G/1} \mathbb{E}[\servicereq_{M/G/1}] =  \arrivalrate \mathbb{E} [\min \{ \servicereq_{1},\dots,\servicereq_{\nrreplicas}\}] < 1. \QEDM
\end{align*}

In case $\nrservers = \nrreplicas$, the above condition is not only sufficient but in fact also necessary since the system behaves exactly as the corresponding $M/G/1$ queue, see also \cite{PC-CERPS}.
In case $\nrservers > \nrreplicas$ the above condition is no longer strictly necessary. However, we will show that, surprisingly, it is asymptotically nearly necessary for independent scaled Bernoulli service requirements, which are defined as 
\begin{align*}
\servicereq = 
\begin{cases}
X \cdot K, \quad & \text{w.p. } 1-p, \\
0, \quad & \text{w.p. } p,
\end{cases}
\end{align*}
where $K$ is a fixed positive real number, and $X$ is a general strictly positive random variable with $\mathbb{E}[X] = 1$. 
Moreover, we assume that $\mathbb{E}[\servicereq]=1$, which implies that $p = 1 - 1/K$.

For notational convenience we label jobs for which none of the $\nrreplicas$ replicas have service requirement $0$ as type-$A$ jobs. For a type-$A$ job $(X_{1}K,\dots,X_{\nrreplicas}K)$ are the service requirements of the replicas at the $\nrreplicas$ sampled servers, where the random variables $X_{1},\dots,X_{\nrreplicas}$ are i.i.d.\ copies of $X$. Jobs for which at least one replica but at most $\nrreplicas-1$ replicas have service requirement equal to $0$ are called type-$B$ jobs, and jobs for which all $\nrreplicas$ replicas have service requirement equal to $0$ are called type-$C$ jobs.\\ 

From Proposition \ref{propo: stability condition IR model} it follows that for independent scaled Bernoulli service requirements the sufficient stability condition reduces to,
\begin{align}
(1-p)^{\nrreplicas} \arrivalrate \mathbb{E}[ \min \{ X_{1}K,\dots, X_{\nrreplicas}K \}] = \frac{\arrivalrate \mathbb{E} [\min \{ X_{1},\dots,X_{\nrreplicas}\}]}{ K^{\nrreplicas-1}} < 1,
\label{eq: sufficient stability condition}
\end{align}
since all jobs, other than type-$A$ jobs, which have arrival rate $(1-p)^{\nrreplicas} \arrivalrate$ and service requirement $\min \{ X_{1}K,\dots, X_{\nrreplicas}K \}$, have service requirements for which \\$\min \{ \servicereq_{1},\dots, \servicereq_{\nrreplicas} \} = 0$.

\section{Asymptotically necessary stability condition}
\label{sec: stability condition}
In this section we shall prove that the sufficient stability condition \eqref{eq: sufficient stability condition} is in fact also asymptotically nearly necessary. The proof relies on the property that the system is most of the time in synchronicity as $K$ grows large.

\begin{figure}[]
    \centering
    \includegraphics[width=7cm]{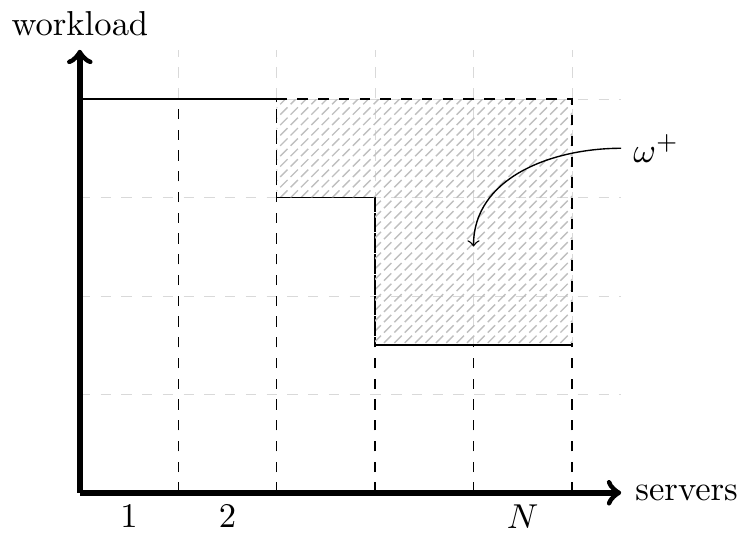} 
    \caption{Visual representation workload surplus} 
    \label{fig: visual surplus workload}  
\end{figure}

In preparation for the proof let us first define a measure for synchronicity. Let the surplus workload, denoted by $\omega^{+}$, be the sum of the (element-wise) differences between the maximum workload and the workload at server $i$ for $i=1,\dots,\nrservers$, i.e., $\omega^{+} = \sum_{i=1}^{\nrservers} \left( \workloadserver_{(1)} - \workloadserver_{i} \right)$; see Figure~\ref{fig: visual surplus workload} for a visual representation. Note that $\omega^{+} = 0$ if and only if the system is in synchronicity. 

In order to prove that the system is in synchronicity nearly all the time, we introduce an auxiliary system which is the same as our system except for three differences. In the auxiliary system (i)~the workload at each server only decreases over time when in synchronicity, (ii)~all type-$A$ jobs are allocated to the first $\nrreplicas$ ordered servers and (iii)~only specific type-$B$ jobs, so-called type-$B_{1}$ jobs, are considered and the other type-$B$ jobs are omitted. We define type-$B_{1}$ jobs as ones for which $\nrreplicas-1$ replicas, with at least one replica with service requirement equal to $0$, are allocated to the first $\nrreplicas-1$ ordered servers and one replica with service requirement $X_{\nrreplicas} K$ to the $\nrservers$-th ordered server, i.e., the server with the lowest current workload.

Below we comment on the properties of the surplus workload $\tilde{\omega}^{+}$ in the auxiliary system.
\begin{property}
\label{prop: omega decrease} The surplus workload in the auxiliary system $\tilde{\omega}^{+}$ experiences downward jumps at the instants of a Poisson process of rate $\frac{(\nrservers-\nrreplicas)!}{\nrservers!} (1-p)p^{\nrreplicas-1} \arrivalrate$, which is exactly the arrival rate of type-$B_{1}$ jobs. The sizes of the downward jumps are equal to $\min \{ \workloadserveraux_{(1)} - \workloadserveraux_{(\nrservers)}, X_{\nrreplicas} K \}$.
\end{property}
Note that the surplus workload in the original system $\omega^{+}$ experiences downward jumps at a higher rate than $\tilde{\omega}^{+}$, since not only type-$B_{1}$ jobs decrease the surplus workload. Moreover, the sizes of the downward jumps in the surplus workload and in the surplus workload in the auxiliary system can differ, since these depend on the workloads in both systems (which are not necessarily equal).

\begin{property}
\label{prop: omega increase} 
The surplus workload in the auxiliary system $\tilde{\omega}^{+}$ experiences upward jumps of size exactly $(\nrservers-\nrreplicas) \min \{ X_{1},\dots,X_{\nrreplicas}\} K$ as a Poisson process of rate $(1-p)^{\nrreplicas}\arrivalrate$, which is the arrival rate of type-$A$ jobs.
\end{property}
Note that the surplus workload in the original system $\omega^{+}$ experiences upward jumps of smaller or equal size, since type-$A$ jobs add at most $\min \{ X_{1},\dots,X_{\nrreplicas}\} K$ work to the current maximum workload, see Remark \ref{rem: workload bound synchronicity}.\\ 

The number of jumps, denoted by $Z$, to reach synchronicity in the auxiliary system when only considering downward jumps is equal to the total number of type-$B_{1}$ jobs that are needed at each server to bridge the difference between the maximum workload and the workload at this server. 
Thus, the expectation of the number of jumps to reach synchronicity, when only considering type-$B_{1}$ jobs and starting in the initial workload state $\workloadstateaux$, where $\workloadstateaux \in \mathcal{S}_{\text{trun}}$, is
\begin{align}
\mathbb{E}[Z] &= \sum_{i=\nrreplicas+1}^{\nrservers} \mathbb{E}[\min\{ n : X_{1}K + \dots + X_{n}K \geq \workloadserveraux_{1} - \workloadserveraux_{i} \} ] \nonumber \\
& \leq \sum_{i=\nrreplicas+1}^{\nrservers} \mathbb{E}[\min\{ n : X_{1}K + \dots + X_{n}K \geq \tilde{\omega}^{+} \} ] \nonumber \\
&\leq (\nrservers-\nrreplicas) \left( \mathbb{E}[\max\{ n : S_{n} \leq \frac{\tilde{\omega}^{+}}{K} \} ] +1 \right) \nonumber \\
&= (\nrservers - \nrreplicas) \left( m(\frac{\tilde{\omega}^{+}}{K}) + 1 \right),
\label{eq: bound number decreases}
\end{align}
where $S_{n} = \sum_{j=1}^{n} X_{j}$ and the renewal function $m$ (cf.~\cite[Def. 10.1.6]{GS-PRP}) is given by $m(t) = \mathbb{E}[ N(t)]$ with $N(t) = \max\{n: S_{n} \leq t \}$. Note that the third line holds with equality if $\frac{\tilde{\omega}^{+}}{K} \not\in \mathbb{N}$.

For proving an asymptotically necessary stability condition, we first need to prove the following two lemmas.
Lemma \ref{lem: sample path wise domination} states that the surplus workload in the auxiliary system stochastically dominates the surplus workload in the original system and Lemma \ref{lem: perfect synchronicity} states that the surplus workload in the auxiliary system is a high fraction of the time equal to $0$ in the long term as $K$ grows large. Together Lemmas \ref{lem: sample path wise domination} and \ref{lem: perfect synchronicity} imply that the original system will also be in synchronicity a high fraction of the time in the long term as $K$ grows large. This in turn implies that almost every arriving job will add $\servicereq_{M/G/1} = \min\{\servicereq_{1},\dots,\servicereq_{\nrreplicas} \}$ to the maximum workload. Observe that this is exactly the upper bound, see Lemma \ref{lem: bound maximum workload}, which resulted in the sufficient stability condition.

Let $\{ \workloadstate(t) , \omega^{+}(t) \}_{t \geq 0}$ denote the stochastic process that describes the evolution of the workload vector $\workloadstate = (\workloadserver_{1},\dots,\workloadserver_{\nrservers})$ and the surplus workload $\omega^{+}$ over time.
We further introduce a stochastic process $\{ \tilde{\workloadstate}(t) , \tilde{\omega}^{+}(t) \}_{t \geq 0}$ that describes the evolution of the workload vector $\workloadstateaux = (\workloadserveraux_{1},\dots,\workloadserveraux_{\nrservers})$ and the surplus workload $\tilde{\omega}^{+}$ of the auxiliary system over time with $\tilde{\omega}^{+}(t) = \sum_{i=1}^{\nrservers} \left( \workloadserveraux_{(1)}(t) - \workloadserveraux_{i}(t) \right)$.

\begin{lemma}
\label{lem: sample path wise domination}
The workload vectors of the auxiliary system and the original system satisfy the inequality $\workloadserveraux_{(1)}(t)-\workloadserveraux_{(i)}(t) \geq \workloadserver_{(1)}(t) - \workloadserver_{(i)}(t)$ for all $t \geq 0$ and $i=1,\dots \nrservers$, when both systems experience the same arrivals, the same generic service requirements and start in the same initial workload state $\workloadstateaux \in \mathcal{S}_{\text{trun}}$, i.e., $\workloadstateaux(0)=\workloadstate(0)=\workloadstateaux$ and $\workloadserveraux_{(1)}=\dots=\workloadserveraux_{(\nrreplicas)}$.
\end{lemma}

\noindent \textbf{Proof:}
Since both systems start in the same initial workload state it follows that $\workloadserveraux_{(1)}(0)-\workloadserveraux_{(i)}(0) = \workloadserver_{(1)}(0) - \workloadserver_{(i)}(0)$, for $i=1,\dots,\nrservers$.
Moreover, by Property~\ref{prop: synchronicity first d servers}, it follows that $\tilde{\omega}_{(1)}(t)-\tilde{\omega}_{(i)}(t) = \workloadserver_{(1)}(t) - \workloadserver_{(i)}(t) = 0$ for $t \geq 0$ and $i=1,\dots, \nrreplicas$.
We prove the statement for $i=\nrreplicas+1,\dots,\nrservers$ by induction in time. Assume that $\workloadserveraux_{(1)}(t_{1})-\workloadserveraux_{(i)}(t_{1}) \geq \workloadserver_{(1)}(t_{1}) - \workloadserver_{(i)}(t_{1})$, then it should hold that $\tilde{\omega}_{(1)}(t_{2})-\workloadserveraux_{(i)}(t_{2}) \geq \workloadserver_{(1)}(t_{2}) - \workloadserver_{(i)}(t_{2})$ for $t_{2} > t_{1}$, when considering all the events that can occur between times $t_{1}$ and $t_{2}$.
\begin{itemize}
\item When no arrivals occur, by definition of both systems, only the value of $\omega^{+}(t)$ can decrease over time. Thus, it follows that $\workloadserver_{(1)}(t_{1}) - \workloadserver_{(i)}(t_{1}) \geq \workloadserver_{(1)}(t_{2}) - \workloadserver_{(i)}(t_{2})$ (which is a strict inequality in case of $\workloadserver_{(i)}(t_{2})=0$), whereas $\workloadserveraux_{(1)}(t_{1})-\workloadserveraux_{(i)}(t_{1}) = \workloadserveraux_{(1)}(t_{2})-\workloadserveraux_{(i)}(t_{2})$. 

\item In case of an arrival of a type-$A$ job the value of $\tilde{\omega}^{+}(t)$ increases with exactly $(\nrservers-\nrreplicas) \min \{ X_{1},\dots,X_{\nrreplicas}\} K$, whereas the value of $\omega^{+}(t)$ increases with at most $(\nrservers-\nrreplicas) \min \{ X_{1},\dots,X_{\nrreplicas}\} K$, see the proof of Lemma~\ref{lem: bound maximum workload} and Property~\ref{prop: omega increase}. Also, note that a type-$A$ job in the auxiliary system is always allocated to the first $\nrreplicas$ ordered servers, instead of $\nrreplicas$ servers sampled uniformly at random. Thus, it follows that $\min \{ X_{1},\dots,X_{\nrreplicas}\} K = \workloadserveraux_{(1)}(t_{2}) - \workloadserveraux_{(1)}(t_{1}) \geq \workloadserver_{(1)}(t_{2}) - \workloadserver_{(1)}(t_{1})$ and $0 =\workloadserveraux_{(i)}(t_{2}) - \workloadserveraux_{(i)}(t_{1}) \leq \workloadserver_{(i)}(t_{2}) - \workloadserver_{(i)}(t_{1})$. Combining the latter two inequalities yields  $\workloadserveraux_{(1)}(t_{2})-\workloadserveraux_{(i)}(t_{2}) \geq \workloadserver_{(1)}(t_{2})-\workloadserver_{(i)}(t_{2})$.

\item In case of an arrival of a type-$B$ job, excluding a type-$B_{1}$ job, only the value of $\omega^{+}(t)$ can decrease. Thus, it follows that $\workloadserver_{(1)}(t_{1}) - \workloadserver_{(i)}(t_{1}) \geq \workloadserver_{(1)}(t_{2}) - \workloadserver_{(i)}(t_{2})$ (which is a strict inequality in case of a type-$B$ job that adds workload to server $i$), whereas $\workloadserveraux_{(1)}(t_{1})-\workloadserveraux_{(i)}(t_{1}) = \workloadserveraux_{(1)}(t_{2})-\workloadserveraux_{(i)}(t_{2})$. 

\item In case of an arrival of a type-$B_{1}$ job the value of $\tilde{\omega}^{+}(t)$ decreases with $\min \{ \workloadserveraux_{(1)}(t_{1}) - \workloadserveraux_{(\nrservers)}(t_{1}), X_{\nrreplicas} K \}$, whereas the value of $\omega^{+}(t)$ decreases with $\min \{ \workloadserver_{(1)}(t_{1}) -  \workloadserver_{(\nrservers)}(t_{1}), X_{\nrreplicas} K \}$. Observe that the decrement in the value of $\tilde{\omega}^{+}(t)$ can be greater than the decrement in the value of $\omega^{+}(t)$, see Property \ref{prop: omega decrease}, but only if $\workloadserver_{(1)}(t_{2}) - \workloadserver_{\nrservers*}(t_{2}) = 0$, where $\nrservers*$ is the server at time $t_{2}$ that had the minimum workload at time $t_{1}$ (which is not necessarily the server with minimum workload at time $t_{2}$). Therefore, it follows that $\workloadserveraux_{(1)}(t_{2})-\workloadserveraux_{(i)}(t_{2}) \geq \workloadserver_{(1)}(t_{2})-\workloadserver_{(i)}(t_{2})$. 
\end{itemize} 
We conclude that in all scenarios it still holds that $\workloadserveraux_{(1)}(t_{2})-\workloadserveraux_{(i)}(t_{2}) \geq \workloadserver_{(1)}(t_{2}) - \workloadserver_{(i)}(t_{2})$. \QED\\

\begin{figure}[]
    \centering
    \includegraphics[width=8.5cm]{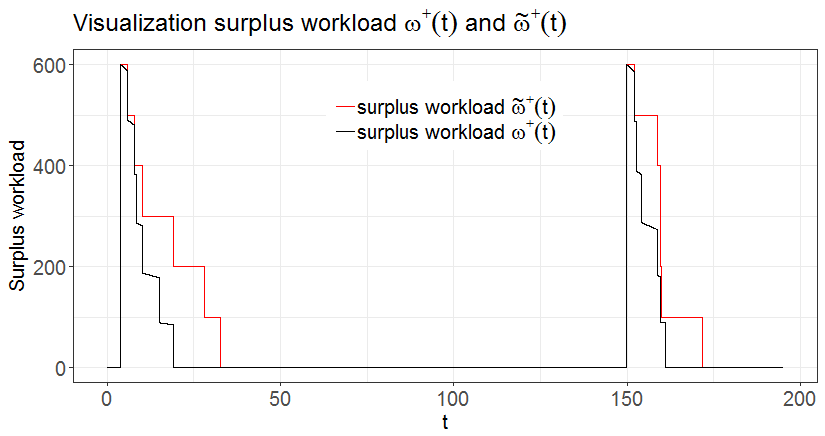} 
    \caption{Visual representation surplus processes, where $\arrivalrate=50$, $\nrservers=8$, $\nrreplicas=2$ and $K=100$.} 
    \label{fig: visual surplus process}  
\end{figure}

Lemma~\ref{lem: sample path wise domination} implies that the surplus workload in the auxiliary system $\tilde{\omega}^{+}(t)$ stochastically dominates the surplus workload $\omega^{+}(t)$ when starting in the same initial workload state, see Figure \ref{fig: visual surplus process}. Now we prove that the surplus workload in the auxiliary system, $\tilde{\omega}^{+}(t)$, is a high fraction of the time equal to $0$ in the long term as $K$ grows large.

\begin{lemma}
\label{lem: perfect synchronicity}
For every $\epsilon>0$ there exists a $K_{\epsilon}(\nrreplicas,\nrservers)$ such that for all $K > K_{\epsilon}(\nrreplicas,\nrservers)$ the value of $ \tilde{\omega}^{+}(t)$ is at least a fraction $(1-\epsilon)$ of the time equal to $0$ in the long term.
\end{lemma}

\noindent \textbf{Proof:} 
First denote $\tau_{1} :=\inf \{t \geq 0 |  \tilde{\omega}^{+}(t)>0 \}$ as the time that the value of $\tilde{\omega}^{+}(t)$ remains equal to $0$, when starting in synchronicity. Note that $\tau_{1}$ is the time until the next upward jump, see Property~\ref{prop: omega increase}. Therefore the expectation of $\tau_{1}$ is given by
\begin{align*}
\mathbb{E}[\tau_{1}] = \frac{1}{(1-p)^{
\nrreplicas}\arrivalrate} = \frac{K^{\nrreplicas}}{\arrivalrate}.
\end{align*}
Denote the time that the workload in the auxiliary system remains in non-synchronicity, i.e., the time that $\tilde{\omega}^{+}(t) > 0$ when starting in initial workload state $\workloadstateaux(0) = \workloadstateaux$, where $\workloadstateaux \in \mathcal{S}_{\text{trun}}$, by $\tau_{2} := \inf \{t \geq 0 |  \tilde{\omega}^{+}(t)=0 \}$. 
Moreover, let $\{ Y | \workloadstateaux(0) = \workloadstateaux \}$ denote the number of increments in the value of $\tilde{\omega}^{+}(t)$ before reaching synchronicity when starting in $\workloadstateaux \in \mathcal{S}_{\text{trun}}$, then the expectation of $\tau_{2}$ is
\begin{align*}
\mathbb{E}[\tau_{2}] &= \sum_{n=0}^{\infty} \mathbb{E}[\tau_{2} | Y=n, \workloadstateaux(0) = \workloadstateaux ] \cdot \mathbb{P}( Y=n | \workloadstateaux(0) = \workloadstateaux ) \\
&= \sum_{n=0}^{\infty} \frac{\mathbb{E}[Z | Y=n, \workloadstateaux(0) = \workloadstateaux ]}{\frac{(\nrservers-\nrreplicas)!}{\nrservers!}(1-p)p^{\nrreplicas-1}\arrivalrate} \cdot \mathbb{P}( Y=n | \workloadstateaux(0) = \workloadstateaux ) \\
& \leq \frac{1}{\frac{(\nrservers-\nrreplicas)!}{\nrservers!}(1-p)p^{\nrreplicas-1}\arrivalrate} \bigg[ (\nrservers-\nrreplicas) \big( m(\frac{\tilde{\omega}^{+}}{K})+1 \big) \\
& \qquad +  \sum_{n=1}^{\infty} \Big( n (\nrservers-\nrreplicas) \big( m(\mathbb{E}[\min \{ X_{1},\dots,X_{\nrreplicas}\}]) + 1 \big) \cdot \mathbb{P}( Y=n | \workloadstateaux(0) = \workloadstateaux ) \Big) \bigg] \\
&= \frac{1}{\frac{(\nrservers-\nrreplicas)!}{\nrservers!}(1-p)p^{\nrreplicas-1}\arrivalrate} \Big[  (\nrservers-\nrreplicas) \big( m(\frac{\tilde{\omega}^{+}}{K})+1 \big) \\
&\quad + (\nrservers-\nrreplicas) \big( m(\mathbb{E}[\min \{ X_{1},\dots,X_{\nrreplicas}\}])+1 \big) \mathbb{E}[Y | \workloadstateaux(0) = \workloadstateaux ] \Big],
\end{align*}
with $\mathbb{E}[\min \{ X_{1},\dots,X_{\nrreplicas}\}] \leq \mathbb{E}[X] = 1$.
The second equality results from Wald's equation, i.e., the equality between the expected time to reach synchronicity (given the number of upward jumps) and the expected number of downward jumps (given the number of upward jumps) multiplied with the expected time between such downward jumps. 
The inequality in the next step results from the proof of Lemma~\ref{lem: bound maximum workload}, which implies that the surplus workload increases with at most $(\nrservers-\nrreplicas) \min \{ X_{1},\dots,X_{\nrreplicas}\} K$ per upward jump, and using the bound on the expected number of downward jumps (given the number of upward jumps), i.e., Equation~\eqref{eq: bound number decreases}.

Together with Wald's equation
\begin{align*}
\mathbb{E}[ Y | \workloadstateaux(0) = \workloadstateaux ] = \mathbb{E}[\tau_{2} ] (1-p)^{\nrreplicas} \arrivalrate,
\end{align*}
we can bound the expected time in non-synchronicity, namely
\begin{align*}
& \mathbb{E}[\tau_{2}] \frac{(\nrservers-\nrreplicas)!}{\nrservers!}(1-p)p^{\nrreplicas-1}\arrivalrate \\ &
\leq (\nrservers-\nrreplicas) \big( m(\frac{\tilde{\omega}^{+}}{K}) + 1 \big) + (\nrservers-\nrreplicas) \big( m(\mathbb{E}[\min \{ X_{1},\dots,X_{\nrreplicas}\}] ) +1 \big) \mathbb{E}[\tau_{2}] (1-p)^{\nrreplicas} \arrivalrate \\
 \Leftrightarrow \quad &\mathbb{E}[\tau_{2}]\left( \frac{(\nrservers-\nrreplicas)!}{\nrservers!}(1-p)p^{\nrreplicas-1}\arrivalrate - (\nrservers-\nrreplicas) \big( m(\mathbb{E}[\min \{ X_{1},\dots,X_{\nrreplicas}\}] ) + 1 \big) (1-p)^{\nrreplicas} \arrivalrate \right) \\
& \leq (\nrservers-\nrreplicas) \big( m(\frac{\tilde{\omega}^{+}}{K}) +1 \big) \\
\Leftrightarrow \quad &\mathbb{E}[\tau_{2}] \leq \frac{ (\nrservers-\nrreplicas) \big( m(\frac{\tilde{\omega}^{+}}{K}) + 1 \big) }{ \frac{(\nrservers-\nrreplicas)!}{\nrservers!}\frac{\arrivalrate}{K}(1-\frac{1}{K})^{\nrreplicas-1} - (\nrservers-\nrreplicas) \big( m(\mathbb{E}[\min \{ X_{1},\dots,X_{\nrreplicas}\}] ) + 1 \big) \frac{\arrivalrate}{K^{\nrreplicas}}} = \frac{1}{\arrivalrate} \mathcal{O}(K),
\end{align*}
under the assumption that $\frac{(\nrservers-\nrreplicas-1)!}{\nrservers!}(1-\frac{1}{K})^{d-1} > \frac{ m(\mathbb{E}[\min \{ X_{1},\dots,X_{\nrreplicas}\}])+1}{K^{\nrreplicas-1}}$. 
Moreover $m(\frac{\tilde{\omega}^{+}}{K}) \downarrow 0$ as $K$ grows large and by renewal theory (cf. \cite{GS-PRP}) we know that \\$m(\mathbb{E}[\min \{ X_{1},\dots,X_{\nrreplicas}\}] ) < \infty$ since $\mathbb{E}[\min \{ X_{1},\dots,X_{\nrreplicas}\}] \leq 1$.

Now if we choose $K_{\epsilon}(\nrreplicas,\nrservers)$ such that for $K=K_{\epsilon}(\nrreplicas,\nrservers)$ one has
\begin{align*}
\frac{\mathbb{E}[\tau_{2}]}{\mathbb{E}[\tau_{1}]+\mathbb{E}[\tau_{2}]} \leq \epsilon,
\end{align*}
then for all $K > K_{\epsilon}(\nrreplicas,\nrservers)$ it follows that 
\begin{align*}
\frac{\mathbb{E}[\tau_{1}]}{\mathbb{E}[\tau_{1}]+\mathbb{E}[\tau_{2}]} > 1 - \epsilon = 1 - \mathcal{O}(\frac{1}{K^{\nrreplicas-1}}).
\end{align*}
This completes the proof that the auxiliary surplus workload is at least a fraction $(1-\epsilon)$ of the time equal to $0$ in the long term.\QED

\begin{remark}
So far it has been assumed that in the initial state the first $\nrreplicas$ ordered workloads are equal, but this assumption is not necessary. One can show via an approach analogous to Lemma~\ref{lem: perfect synchronicity}, but with bound $\mathbb{E}[Z ] < \nrservers \big( m(\frac{\tilde{\omega}^{+}}{K}) + 1 \big)$, that the expected time to reach synchronicity when starting in an arbitrary initial workload state is still finite. Note that after reaching synchronicity the assumption is valid and that directly after synchronicity $\tilde{\omega}^{+}(t)=\omega^{+}(t) = (\nrservers-\nrreplicas) \min \{ X_{1},\dots,X_{\nrreplicas}\}K$.
\end{remark}

Now we are ready to prove the main theorem of the paper.

\begin{theorem}
\label{thm: necessary stability condition}
For every $\epsilon>0$ there exists a $K_{\epsilon}(\nrreplicas,\nrservers)$ such that for all $K > K_{\epsilon}(\nrreplicas,\nrservers)$ a necessary stability condition for independent scaled Bernoulli service requirements is
\begin{align}
(1-\epsilon) \frac{\arrivalrate \mathbb{E}[ \min \{ X_{1},\dots,X_{\nrreplicas}\}]}{ K^{d-1}} < 1.
\label{eq: non stable necessary condition}
\end{align}
\end{theorem}
\noindent \textbf{Proof:} From Lemma~\ref{lem: sample path wise domination} we know that $\tilde{\omega}^{+}(t)$ stochastically dominates $\omega^{+}(t)$ and Lemma~\ref{lem: perfect synchronicity} states that for every $\epsilon>0$ there exists a $K_{\epsilon}(\nrreplicas,\nrservers)$ such that for all $K > K_{\epsilon}(\nrreplicas,\nrservers)$ the value of $\tilde{\omega}^{+}(t)$ is at least a fraction $(1-\epsilon)$ of the time equal to $0$ in the long term. 
Hence this latter statement also holds for the value of $\omega^{+}(t)$. 
Moreover, by definition, if $\omega^{+}(t) = 0$ then the system is in synchronicity.
In synchronicity, type-$A$ jobs add exactly $\min \{ X_{1},\dots,X_{\nrreplicas}\} K$ work to the sampled servers. 
We conclude that, independent of the behavior in non-synchronicity, in the long term at least a fraction $(1-\epsilon)$ of the type-$A$ jobs adds exactly $\min \{ X_{1},\dots,X_{\nrreplicas}\} K$ work to the current maximum workload. Thus, for the system to be stable it should at least be able to handle these latter type-$A$ jobs. 
\QED

\begin{remark}
The expected time in non-synchronicity depends on the renewal function $m(t)$, see Lemma~\ref{lem: perfect synchronicity}. This function in turn depends on the distribution of the $X$ component in the service requirement distribution $\servicereq$. For some distributions an explicit expression for $m(t)$ is known (cf.~\cite{GS-PRP}):
\begin{itemize}
\item $X \equiv 1$:
$m(t) = \lfloor t \rfloor$,
\item $X \sim \text{Exp}(1)$: 
$m(t) = t$,
\item $X \sim \text{Unif}[0,2]$:
$m(t)+1 = \sum_{i=0}^{\lfloor t/2 \rfloor} (-1)^{i} \frac{(t/2-i)^{i}}{i!} e^{t/2-i}$.
\end{itemize}
\end{remark}

\section{Numerical results}
\label{sec: numerical results finite K}
In Section~\ref{sec: stability condition} it is proven that the system, for scaled Bernoulli distributed service requirements and $K$ large enough, is a high fraction of the time in synchronicity in the long term. In this section we will use simulation to quantify this statement for various values of $\nrservers$ and $K$, where $\nrreplicas=2$ is fixed.
  
\begin{figure}[]
\centering
\includegraphics[width=8cm]{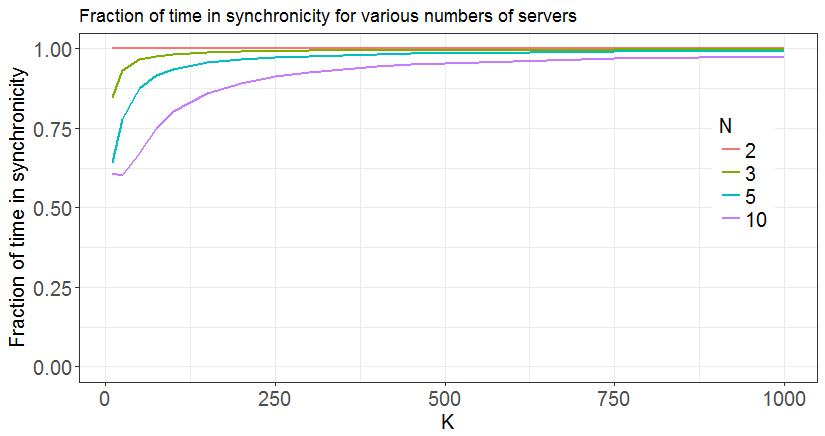}
\includegraphics[width=8cm]{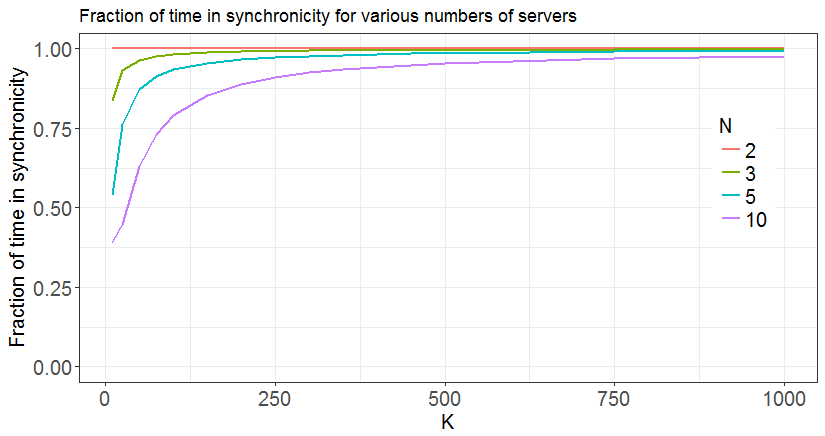}
\caption{Long-term fraction of time in synchronicity (obtained by simulation) for the setting $\nrreplicas=2$, $X \equiv 1$ and $\frac{\arrivalrate(K)}{K} = 0.5$ (top) and $\frac{\arrivalrate(K)}{K} = 0.9$ (bottom).}
\label{fig: simulated synchronicity prob}
\end{figure}

In Figure~\ref{fig: simulated synchronicity prob} the long-term fraction of time in synchronicity is depicted as a function of $K$ for various values of $\nrservers$, where we allow $\arrivalrate$ to depend on $K$ and write $\arrivalrate(K)$ to reflect that. 
It can be seen that the system with $\nrservers=\nrreplicas$ is always in synchronicity, which follows from Property~\ref{prop: synchronicity first d servers}. 
Moreover, the long-term fraction of time in synchronicity is higher for lower values of $\frac{\arrivalrate(K)}{K}$. The reason is that the empty state is included in the definition of synchronicity. 
Another observation is that for fixed $\arrivalrate(K)$ and $K$, increasing $\nrservers$ decreases the long-term fraction of time in synchronicity. 
This is related to the fact that $K_{\epsilon}(\nrreplicas, \nrservers)$ defined in Theorem~\ref{thm: necessary stability condition} depends on $\nrservers$.\\

In Lemma~\ref{lem: bound maximum workload} we proved that the maximum workload is bounded by the workload in a corresponding $M/G/1$ queue. 
From this lemma it follows that for independent scaled Bernoulli service requirements, the maximum workload is bounded by the workload in a corresponding $M/G/1$ queue with arrival rate $\arrivalrate_{M/G/1}(K) = (1-p)^{\nrreplicas}\arrivalrate$ and service requirement $\servicereq_{M/G/1}(K) =\min \{ X_{1},\dots,X_{\nrreplicas}\}K$ since all arrivals, other than the arrivals of type-$A$ jobs, have service requirements for which $\min \{ \servicereq_{1},\dots, \servicereq_{\nrreplicas} \} = 0$.
This bound can be used to find an upper bound on the expected waiting time since an arriving job needs to wait at most for the current maximum workload, which is bounded by the workload $\workload_{M/G/1}$ in the corresponding $M/G/1$ queue. From $M/G/1$ theory (cf.~\cite[Section X.3]{A-APQ}) we get
\begin{align}
\mathbb{E}[\waitingtime] &\leq \mathbb{E}[\workload_{M/G/1}] \nonumber \\
& = \frac{\arrivalrate_{M/G/1}(K) \mathbb{E}[\servicereq_{M/G/1}^2(K)]}{2 ( 1- \arrivalrate_{M/G/1}(K) \mathbb{E}[\servicereq_{M/G/1}(K)])} \nonumber \\
&= \frac{(1-p)^{\nrreplicas}\arrivalrate \mathbb{E}[\min \{ X_{1},\dots,X_{\nrreplicas} \}^2] K^2}{2 ( 1- (1-p)^{\nrreplicas}\arrivalrate \mathbb{E}[\min \{ X_{1},\dots,X_{\nrreplicas} \}]  K)}.
\label{eq: upper bound EW}
\end{align}
Note that this bound is tight for $ \nrservers = \nrreplicas$ since the system behaves exactly as the corresponding $M/G/1$ queue and asymptotically tight in $K$ for $\nrservers > \nrreplicas$.
For $X \equiv 1$ constant and $\nrreplicas=2$ we get
\begin{align*}
\mathbb{E}[\waitingtime] &\leq \frac{(1-p)^{2}\arrivalrate K^2}{2 ( 1- (1-p)^{2}\arrivalrate K)}\\
&= \frac{\arrivalrate}{2 (1 - \frac{\arrivalrate}{K})},
\end{align*}
which is linear in $K$ if we assume that $ \frac{\arrivalrate}{K}$ is fixed. Notice that this upper bound does not depend on the number of servers.

\begin{figure}[]
\centering
\includegraphics[width=8cm]{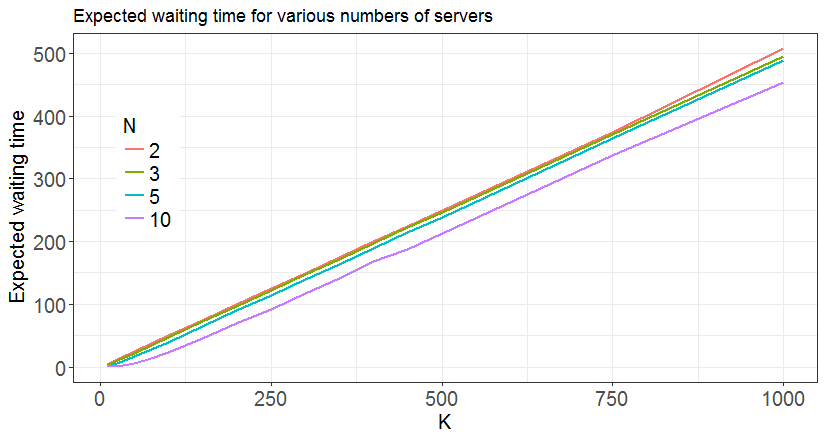}
\includegraphics[width=8cm]{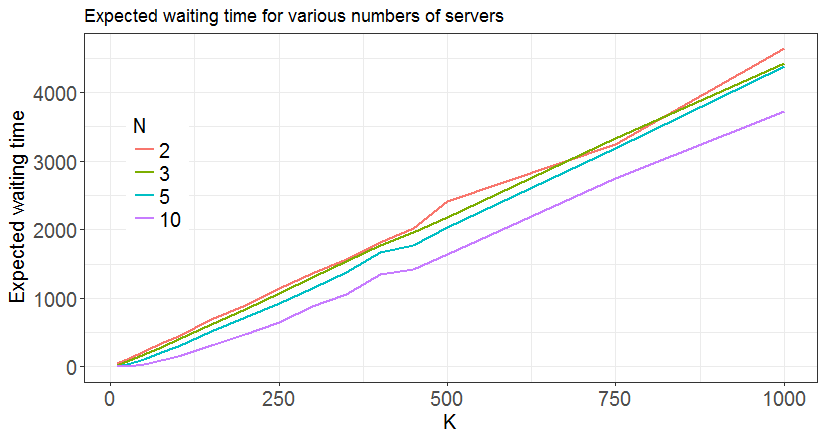}
\caption{Expected waiting time (obtained by simulation) for the setting $\nrreplicas=2$, $X \equiv 1$ and $\frac{\arrivalrate(K)}{K} = 0.5$ (top) and $\frac{\arrivalrate(K)}{K} = 0.9$ (bottom). Note that $\nrservers=2$ corresponds to the upper bound given in \eqref{eq: upper bound EW}.}
\label{fig: simulated expected waiting time}
\end{figure}

Figure~\ref{fig: simulated expected waiting time} shows  the expected waiting time as a function of $K$ for various values of $\nrservers$; again we allow $\arrivalrate$ to depend on $K$ and write $\arrivalrate(K)$ to reflect that. When comparing both figures we can conclude that for finite $K$ the number of servers $\nrservers$ influences the expected waiting time more than the value of the fraction $\frac{\arrivalrate(K)}{K}$. Moreover, for $\nrservers$ large, also a larger $K$ is needed for the upper bound to be accurate. 

Observe that with the upper bound for the expected waiting time we also have an upper bound for the expected latency, since
\begin{align*}
\mathbb{E}[\latency] &= \mathbb{E}[\waitingtime] + \mathbb{E}[\min \{ \servicereq_{1},\dots,\servicereq_{\nrreplicas}\}],
\end{align*}
where $\mathbb{E}[\min \{ \servicereq_{1},\dots,\servicereq_{\nrreplicas}\}] \leq 1$ since by assumption $\mathbb{E}[B_{i}] = 1$, for $i=1,\dots,\nrreplicas$.

\begin{figure}[]
\centering
\includegraphics[width=8cm]{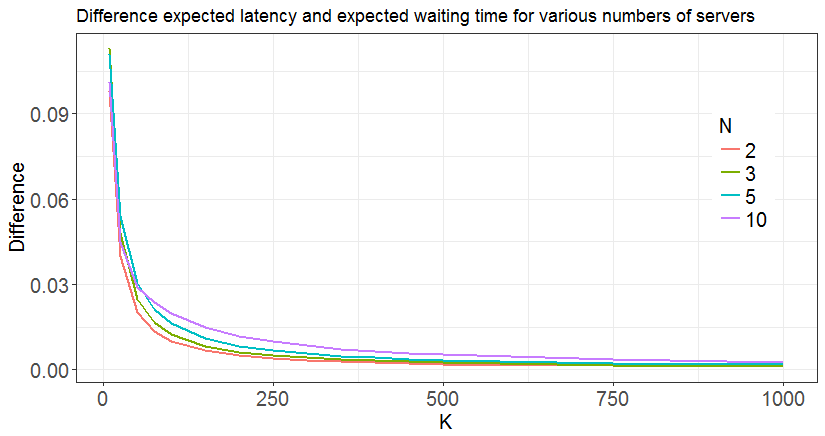}
\caption{Difference between the expected latency and the expected waiting time (obtained by simulation) for the setting $\nrreplicas=2$, $X \equiv 1$ and $\frac{\arrivalrate(K)}{K} = 0.5$.}
\label{fig: simulated difference expected latency waiting time}
\end{figure}

Note that the upper bound for the service requirements, i.e.,  $\mathbb{E}[\min \{ \servicereq_{1},\dots,\servicereq_{\nrreplicas}\}] \leq 1$, is not (asymptotically) tight in $K$, since $\mathbb{E}[\min \{ \servicereq_{1},\dots,\servicereq_{\nrreplicas}\}] \downarrow 0$ as $K$ grows large. In Figure \ref{fig: simulated difference expected latency waiting time} the difference between the expected latency and the expected waiting time is depicted as function of $K$ for various values of $\nrservers$. Indeed, it can be seen that this difference vanishes as $K$ grows large.

\section{Conclusion}
\label{sec: conclusion}
In this paper we have proven that the maximum workload in a parallel-server system with c.o.c.\ redundancy is upper bounded by the workload in a related $M/G/1$ queue. This directly yields a sufficient stability condition. Moreover, we proved that in the case of independent scaled Bernoulli service requirements the system is asymptotically a high fraction of the time in so-called synchronicity in the long term. In synchronicity the upper bound of the related $M/G/1$ queue is in fact tight, and this resulted in an asymptotically necessary stability condition. 
Interestingly, both the sufficient and asymptotically nearly necessary condition is independent of the number of servers, but do depend on the number of replicas $\nrreplicas$. 
In contrast, in the case of exponentially distributed service requirements the stability condition depends linearly on the number of servers and not on the number of replicas. This indicates that the stability condition in a c.o.c.\ redundancy system with i.i.d.\ service requirements is highly sensitive to the distribution of these service requirements. 

The bound on the maximum workload also resulted in an upper bound for the expected waiting time, which is again (asymptotically) tight (as the scale of the service requirement grows large). This bound directly resulted in an upper bound for the expected latency. 

We assumed that jobs arrive according to a Poisson process, but it might be possible to relax this assumption. In particular, the proof of the sufficient stability condition does not rely on Poisson arrivals and could be extended to a general arrival process. The extension of the proof of the asymptotically necessary condition is more involved. Another interesting topic for further research is to extend the developed framework to obtain the stability condition for more general service requirements. 




\end{document}